\documentclass[11pt, reqno]{amsart}
\usepackage[english]{babel}
\usepackage[table]{xcolor}
\usepackage[utf8]{inputenc}
\usepackage{amsmath, amsfonts, amsthm, amscd, amssymb, fullpage, fancyhdr, upgreek, amssymb, hyperref, enumerate, comment, siunitx, enumitem, graphicx, mathrsfs, mathtools, filecontents, xcolor,
microtype}
\usepackage{listings}
\usepackage[scr]{rsfso}
\usepackage{thmtools}
\usepackage{thm-restate}

\usepackage{tikz}
\usepackage{tkz-tab}
\usepackage{tkz-graph}
\usetikzlibrary{shapes.geometric,positioning}

\numberwithin{equation}{section}

\newcommand{\Z}{\mathbb{Z}}

\newcommand{\eps}{\varepsilon}
\newcommand{\dq}[1]{\lq\lq #1\rq\rq}

\newcommand{\set}[1]{\left\{#1\right\}}
\newcommand{\prob}[1]{\mathbb{P}\left(#1\right)}
\newcommand{\e}[1]{\mathbb{E}\left[#1\right]}
\newcommand{\floor}[1]{\left\lfloor #1\right\rfloor}
\newcommand{\ceil}[1]{\left\lceil #1\right\rceil}

\newcommand{\yt}{\tilde{Y}}
\newcommand{\zt}{\tilde{Z}}
\newcommand{\wt}{\tilde{W}}

\newcounter{dummy} \numberwithin{dummy}{section}
\newtheorem{lemma}[dummy]{Lemma}
\newtheorem{proposition}[dummy]{Proposition}
\newtheorem{theorem}[dummy]{Theorem}
\newtheorem{corollary}[dummy]{Corollary}
\theoremstyle{definition} \newtheorem{definition}[dummy]{Definition}
\theoremstyle{definition}

\newcommand{\hr}[1]{\href{#1}{\url{#1}}}

\title{Limiting Behavior in Missing Sums of Sumsets}

\author{Aditya Jambhale}
\address{University of Cambridge}

\author{Rauan Kaldybayev}
\address{Williams College}

\author{Steven J. Miller}
\address{Williams College}

\author{Chris Yao}
\address{Yale University}


\thanks{This work was completed during the 2023 SMALL REU program at Williams College. It was supported in part by NSF Grants DMS1561945 and DMS1659037, Williams College, and Churchill College, Cambridge.}

\date{\today}

\begin{document}

\maketitle

\begin{abstract}
    We study $|A + A|$ as a random variable, where $A \subseteq \{0, \dots, N\}$ is a random subset such that each $0 \le n \le N$ is included with probability $0 < p < 1$, and where $A + A$ is the set of sums $a + b$ for $a,b$ in $A$. Lazarev, Miller, and O'Bryant studied the distribution of $2N + 1 - |A + A|$, the number of summands not represented in $A + A$ when $p = 1/2$. A recent paper by Chu, King, Luntzlara, Martinez, Miller, Shao, Sun, and Xu generalizes this to all $p\in (0,1)$, calculating the first and second moments of the number of missing summands and establishing exponential upper and lower bounds on the probability of missing exactly $n$ summands, mostly working in the limit of large $N$. We provide exponential bounds on the probability of missing at least $n$ summands, find another expression for the second moment of the number of missing summands, extract its leading-order behavior in the limit of small $p$, and show that the variance grows asymptotically slower than the mean, proving that for small $p$, the number of missing summands is very likely to be near its expected value.
\end{abstract}

\setcounter{tocdepth}{1}
\tableofcontents

\newpage

\section{Introduction}
Fix a real number $0 < p < 1$ and an integer $N \ge 0$, and let $A$ be a random subset of $\set{0, \dots, N}$ such that each $n$ between $0$ and $N$ is included in $A$ independently with probability $p$. We consider the sumset
\begin{equation}
    A + A \ \coloneqq\ \set{a + b: a, b \in A} \subseteq \set{0, \dots, 2N}.
\end{equation}

There is an extensive literature on \dq{more sums than differences} (MSTD) sets, defined as those sets $A$ whose sumset cardinalities $|A + A|$ are greater than difference set cardinalities $|A - A|$, where $A - A$ is the set of elements $a - b$ with $a, b \in A$.

MSTD sets are unusual because addition is commutative while subtraction is anticommutative, and one would expect $A - A$ to usually have more elements than $A + A$; see for example \cite{AMMS, BELM, CKLMMSSX, CLMS, CMMXZ, DKMMW, He, HLM, ILMZ, Ma, MOS, MS, MPR, MV, Na1, Na2, PW, Ru1, Ru2, Ru3, Sp, Zh1}. 

While there were known constructions of infinite families of MSTD sets, these examples had density zero as $N\to\infty$. It was surprising when Martin and O'Bryant \cite{MO} proved that positive fraction of all sets $A \subseteq\{1, \dots, n\}$ are MSTD. Martin and O'Bryant were also interested in the distribution of $|A + A|$ as a random variable and calculated some fundamental results, such as the expectation value of $|A + A|$ in the case $p = 1/2$. Lazarev, Miller, and O'Bryant \cite{LMO} proved further results, such as exponential decay of the probability of missing many summands, and gave a rigorous proof to the observation that sumsets show a noticeable bias towards missing an even number of summands, for $p = 1/2$. The 2020 paper of \cite{CKLMMSSX} generalized \cite{LMO}'s proofs to the case of arbitrary $p$. We cover some of these results while also proving certain new ones. In particular, we carefully study the second moment $\e{Y^2}$ and use it to show that the distribution of missing summands is sharply peaked around the mean for large $N$ and small $p$.

\begin{figure}[h] \label{fig:non-inclusion in A + A}
\caption{Probabilities of non-inclusion of $n$ into $A + A$ for $0 \le n \le 2N$ and $p = 1/2$, $N = 40$.}
\centering
\includegraphics[width=0.6\textwidth]{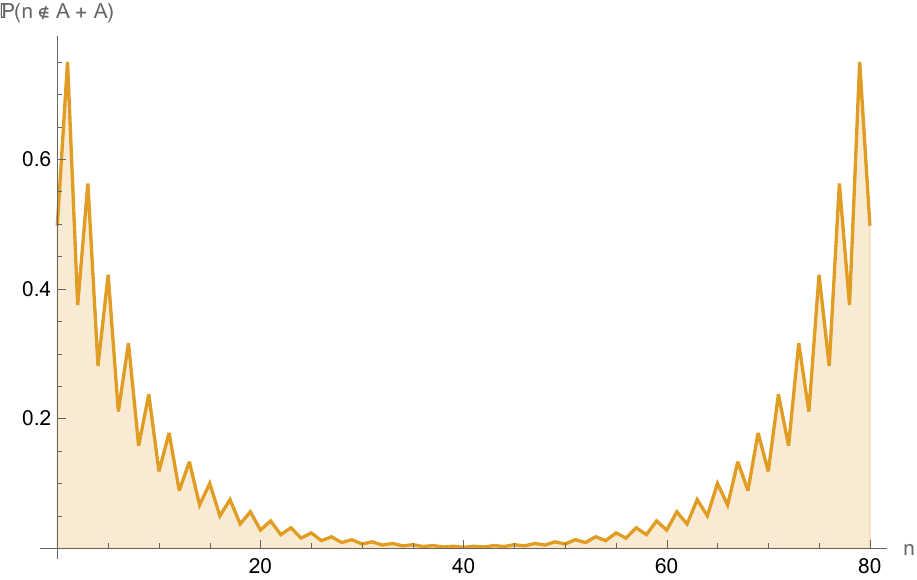}
\end{figure}

We note the following facts about the distribution of $A+A$.
\begin{itemize}
    \item[1.] The set $A + A$ is \dq{almost full} in the middle. While numbers near $N$ are almost certainly included in $A + A$, numbers close to $0$ (the \emph{left fringe}) or $2N$ (the \emph{right fringe}) have a significant chance of being left out \cite{Zh2}. See Figure \ref{fig:non-inclusion in A + A}.
    \item[2.] The set $A + A$ has \dq{symmetry} around $N$ in the sense that $|A+A|$ and $|2N - (A+A)|$ are identically distributed.
    \item[3.] If $j-i > N$, the inclusions of $i$ and $j$ in $A+A$ are independent, but for general $i$ and $j$ they are not.
\end{itemize}

We continue to investigate the number of elements missing from $A + A$. We define the random variables $Y, Z, W$ to count missing summands in the left fringe, right fringe, and the whole of $A + A$ respectively. Let $Y$ be the number of integers at most $ N$ that are missing from $A + A$,
\begin{align} \begin{split}
    Y &\ \coloneqq\ |\set{0 \le n \le N: n \notin A + A}|\\
        &\quad =\ N + 1 - |(A + A) \cap \set{0, \dots, N}|,
\end{split} \end{align}
let $Z$ be the number of integers from $N + 1$ to $ 2N$ that are missing from $A + A$,
\begin{align} \begin{split}
    Z &\ \coloneqq\ |\set{N + 1 \ \le\ n \ \le\ 2 N: n \notin A + A}|\\
        &\quad =\ N - |(A + A) \cap \set{N + 1, \dots, 2 N}|,
\end{split} \end{align}
and $W$ the number of integers up to $2N$ that are missing from $A + A$,
\begin{align} \begin{split}
    W &\ \coloneqq\ |\set{0 \le n \ \le\ 2N: n \notin A + A}|\\
        &= 2 N + 1 - |A + A| = Y + Z.
\end{split} \end{align}

Define $X_n$ as the indicator variable for the \emph{non}-inclusion of $n$ in $A + A$:
\begin{equation}
    X_n \ \coloneqq\ [n \notin A + A] \ =\ \begin{cases}
        1 & n \notin A + A\\
        0 & n \in A + A.
    \end{cases}
\end{equation}

Then,
\begin{equation}Y \ =\ \sum_{n=0}^{N}{X_n}, \quad Z \ =\ \sum_{n=N + 1}^{2N}{X_n}, \quad W \ =\ \sum_{n=0}^{2N}{X_n}.\end{equation}

Note $Y$ and $Z$ have largely the same probability distributions due to the symmetry about $N$, except that $Y$ ranges from $0$ to $N + 1$ while $Z$ from $0$ to $N$. We make use of this symmetry to simplify calculations, proving results about $Y$ and using Theorem \ref{thm:W convolution} to convert to similar statements about $W$. 

As in \cite{CKLMMSSX} and \cite{MO}, we have for $0 \le n \le N$,
\begin{equation} \label{eq:prob n small in A + A}
    \prob{n \in A + A} \ =\ \begin{cases}
        1 - (1 - p^2)^{\frac{n+1}{2}} & \text{$n$ odd}\\
        1 - (1 - p)(1 - p^2)^{n/2} & \text{$n$ even.}
    \end{cases}
\end{equation}

We can use symmetry to write an analogous formula for $N < n \le 2N$ as
\begin{equation} \label{eq:prob n large in A + A}
    \prob{n \in A + A} \ =\ \begin{cases}
        1 - (1 - p^2)^{\frac{(2N-n)+1}{2}} & \text{$n$ odd}\\
        1 - (1 - p)(1 - p^2)^{(2N-n)/2} & \text{$n$ even.}
    \end{cases}
\end{equation}

Recall that $Y = \sum_{n=0}^{N}{X_n}$, so by linearity of expectation,
\begin{align} \begin{split} \label{en:e y}
    \e{Y} &\ =\ \sum_{\text{$n$ even}}{(1 - p)(1 - p^2)^{n/2}} + \sum_{\text{$n$ odd}}{(1 - p^2)^{\frac{n+1}{2}}}\\
        &\ =\ \left(\frac{2}{p^2} - \frac{1}{p} - 1\right) - \left(\sqrt{1 - p^2}\right)^N \begin{cases}
            \frac{(2 - p)(1 - p^2)}{p^2} & \text{$N$ even}\\
            \frac{(2 - p - p^2)\sqrt{1 - p^2}}{p^2} & \text{$N$ odd.}
        \end{cases}
\end{split} \end{align}

For the entire sumset, the expected number of missing summands is
\begin{align} \begin{split} \label{en:e w}
    \e{W} \ =\ \left(\frac{4}{p^2} - \frac{2}{p} - 2\right) - \left(\sqrt{1 - p^2}\right)^N \begin{cases}
            \frac{4 - 2p - 2p^2 + p^3}{p^2} & \text{$N$ even}\\
            \frac{(4 - 2p - p^2)\sqrt{1 - p^2}}{p^2} & \text{$N$ odd.}
        \end{cases}
\end{split} \end{align}

In \textsection\ref{sec: independence of fringes}, we study the number of missing summands in the left and right fringes. It is well known that the left and right fringes are independent when their lengths are smaller than $N/2$ \cite{CKLMMSSX}, but larger sized fringes (of length $\geq N/2$) have not been studied extensively. While the left and right fringes are not truly independent in this case, we prove the following theorem, which intuitively says that $Y$ and $Z$ are ``asymptotically independent''.

\begin{restatable*}{theorem}{convolutionTheorem} 
\label{thm:W convolution}
For any $0 \le m \le 2N$,
\begin{equation}\left|\prob{W = m} - \sum_{y=0}^{m}{\prob{Y=y}\prob{Z=m-y}}\right| \ \le\ \frac{8}{p^2} (1 - p^2)^{N/4}.\end{equation}

We say that as $N \to \infty$, $W$ is the convolution of $Y$ with $Z$.
\end{restatable*}

In \textsection\ref{sec: exponential bounds on missing summands}, we study the $k$\textsuperscript{th} moment of the number of missing summands $\e{Y^k}$ and bound the probability of missing at least $n$ summands. We take a different approach than the authors of \cite{Zh2} and \cite{CKLMMSSX}, who study $m_{n;p}(k)\coloneqq \prob{2N+1 - |A+A| = k}$ and $m_p(k)\coloneqq \lim_{n\to\infty} m_{n;p}(k)$, which allows us to obtain a simpler bound. For $\alpha$ defined in \eqref{eq:alpha defn}, we have the following.

\begin{restatable*}{corollary}{probabilityExponentialBoundCorollary}
\label{eq:Y ge n O bound}
For any real number $\eps > 0$, (and again, for $n > 1 / \alpha$),
\begin{align}
    \prob{Y \ge n} &\ =\ O\left(n \left(\sqrt{1 - p^2}\right)^n\right) \nonumber \\
    &\ =\ O\left(\left(\sqrt{1 - p^2} + \eps\right)^n\right). 
\end{align}
\end{restatable*}

We then improve on this exponential bound in \textsection \ref{sec:twosummands} by using a more detailed estimate for the variance. For suitable $\lambda_1$ as defined in \eqref{eq:defn lambda12}, we have the following.

\begin{restatable*}{corollary}{betterProbabilityBoundCorollary} \label{corr:better bound on prob of missing many}
We have
\begin{equation}\prob{Y \ge n} \ =\ O\left((\lambda_1 + \eps)^n\right)\end{equation}
for all $\eps > 0$. For $p = 1/2$, $\lambda_1 = \varphi/2 \approx 0.81$ where $\varphi$ is the golden ratio. More precisely, for $n > 1 / \alpha'$,
\begin{equation}\prob{Y \ge n} \ \le\ \left(\frac{2}{p^2} - \frac{1}{p} - 1\right) e^{-(n-1)(\alpha' - 1/n)} + \frac{2 n \alpha' e^{-n \alpha' + 1}}{\lambda_1},\end{equation}
where $\alpha' \coloneqq \log{(1 / \lambda_1)}$.
\end{restatable*}

We study the behavior of the second moment $\e{Y^2}$ in the limit as $N\to\infty$ in \textsection\ref{sec: second moment}, providing an exact expression for it. By considering this limit, we find a simpler expression of the second moment than that of \cite{CKLMMSSX}, which finds a closed-form expression for the second moment for finite $N$. The authors hope that further considering the ``infinite case'' (i.e. chosing $A\subseteq \mathbb{N}$) may lead to interesting results in the finite case, which is of broader interest. Finally, in \textsection\ref{sec: concentration and asymptotics}, we conclude by finding the leading order term of the second moment $\e{Y^2}$ in the limit of large $N$ and  determining the asymptotic behavior of $\text{Var}(Y)$.

\begin{restatable*}{proposition}{secondMomentLeadingOrderProp}
\label{eq:second moment leading order} 
There is an error term $\delta(p)$ such that
\begin{equation}
    \lim_{N \to \infty}{\e{Y^2}} \ =\ \frac{4}{p^4} + \delta(p),
\end{equation}
where $\lim_{p \to 0}{\delta(p) / p^{-4}} = 0$.
\end{restatable*}

In \textsection\ref{sec: future work}, we conclude by considering future work and questions that arise naturally from the results presented here.


\section{Independence of the Fringes}\label{sec: independence of fringes}
This section proves that missing summands in the left and right fringes are almost independent. Note that if we had defined the left fringe as integers at most $N/2$ missing from $A+A$ and the right fringe as integers between $3N/2$ and $2N$, we would automatically have independence, but this is not true in our case. We show that the distribution of the total number of missing summands is almost the convolution with itself of the distribution of the number of missing summands in one fringe (Theorem \ref{thm:W convolution}). For the rest of this paper we will focus on counting missing summands in the left fringe, which is convenient. Define random variables $\yt$ and $\zt$ as
\begin{align} \begin{split}
    \yt \ \coloneqq\ \sum_{n=0}^{\floor{N/2}}{X_n}, \qquad \zt \ \coloneqq\ \sum_{n=\floor{3N/2}+1}^{2N}{X_n}
\end{split} \end{align}
to count the number of missing summands on the \dq{very left} and \dq{very right,} respectively. These are exactly independent, since $X_i$ and $X_j$ are independent whenever $|i - j| > N$. (If $j - i > N$, then $0 \le i < N < j \le 2N$, and $X_i$ depends on the inclusion of $0, \dots, i$ into $A$ while $X_j$ depends on the inclusion of $j - N, \dots, N$.) By Lemma \ref{lemma:e dy dz} below, $Y$ is equal to $\yt$ with high probability and $Z$ to $\zt$, and therefore, $Y$ and $Z$ are almost independent.

\begin{lemma} \label{lemma:e dy dz}
The expectation values $\e{Y - \yt}$ and $\e{Z - \zt}$ are bounded above by $\frac{2}{p^2} \left(1 - p^2\right)^{N/4}$ and below by $0$.
\end{lemma}
\begin{proof}
By \eqref{eq:prob n small in A + A}, the probability for the non-inclusion of $i$ in $A + A$ is at most $(1 - p^2)^{(i+1)/2}$. Therefore,
\begin{align} \begin{split}
    \e{Y - \yt} \ \le\ \sum_{i=\floor{N/2}+1}^{N}{\left(\sqrt{1 - p^2}\right)^{i+1}} \ <\ \frac{2}{p^2} \left(\sqrt{1 - p^2}\right)^{N/2},
\end{split} \end{align}
and a similar calculation works for $Z - \zt$. The expectation values are nonnegative because each is a sum of probabilities $\e{X_i}$.
\end{proof}
\begin{corollary}
The probabilities of disagreement $\prob{Y \ne \yt}$ and $\prob{Z \ne \zt}$ are both bounded above by $\frac{2}{p^2} \left(1 - p^2\right)^{N/4}$.
\end{corollary}
\begin{proof}
$Y - \yt$ and $Z - \zt$ can only take nonnegative integer values, and Markov's inequality gives us the desired result.
\end{proof}

\subsection{Convolution of Left and Right Fringes}
Recall that the total number of missing summands $W$ is the sum of the number of missing summands $Y$ on the left and $Z$ on the right. Since the latter two are largely independent, the distribution of $W$ is close to the convolution of $Y$ with $Z$. Theorem \ref{thm:W convolution} is similar to Theorems 6.4 and 6.9 of \cite{CKLMMSSX}.

\begin{lemma} \label{lemma:prob ne}
If $S, T$ are random variables that take values in some finite set $A$, and if $x$ is a possible value they can take, then
\begin{equation}\left|\prob{S = x} - \prob{T = x}\right| \ \le\ \prob{S \ne T}.\end{equation}
\end{lemma}
\begin{proof}
Consider the random variable $(S, T)$, which takes on values in the Cartesian product $A \times A$. Then $\prob{S \ne T}$ is the sum of probabilities $\prob{S = y, T = z}$ over all \dq{off-diagonal} pairs $y \ne z$, while
\begin{align} \begin{split}
    \left|\prob{S = x} - \prob{T = x}\right| &\ =\ \left|\sum_{y \in A}{\prob{S = x, T = y}} - \sum_{y \in A}{\prob{S = y, T = x}}\right|\\
                                            &\ =\ \left|\sum_{y \ne x}{\prob{S = x, T = y} - \sum_{y \ne x}{\prob{S = y, T = x}}}\right|\\
                                            &\ \le\ \sum_{y \ne x}{\prob{S = x, T = y}} + \sum_{y \ne x}{\prob{S = y, T = x}}\\
\end{split} \end{align}
is at most the sum over the column $\set{x} \times A$ and the row $A \times \set{x}$ excepting the diagonal point $(x, x)$. Clearly, this is not greater than the sum over all off-diagonal elements.
\end{proof}

\convolutionTheorem

\begin{proof}
Define $\wt \coloneqq \yt + \zt$. Since $\yt$ and $\zt$ are independent,
\begin{equation} \label{eq:wt eq m prob}
    \prob{\wt = m} \ =\ \sum_{y=0}^{m}{\prob{\yt=y}\prob{\zt=m-y}}.
\end{equation}

To translate this equality involving $\yt, \zt, \wt$ into an inequality involving $Y, Z, W$, we use Lemmas \ref{lemma:e dy dz} and \ref{lemma:prob ne}. For the left-hand side,
\begin{equation} \label{eq:prob Wm prob Wtm almost equal}
    \left|\prob{W = m} - \prob{\wt = m}\right| \ \le\ \prob{W \ne \wt} \ \le\ \e{W - \wt} \ \le\ \frac{4}{p^2}(1 - p^2)^{N/4},
\end{equation}
where we have used that $W - \wt$ only takes on nonnegative integer values and applied the linearity of expectation to find $\e{W - \wt}$ as $\e{Y - \yt} + \e{Z - \zt}$. For the right-hand side, we note that
\begin{align}
    \bigg|\prob{\yt=y}&\prob{\zt=m-y} - \prob{Y=y}\prob{Z=m-y} \bigg| \nonumber \\
    &\le\ \left|\prob{\yt=y} - \prob{Y=y}\right| \prob{\zt=m-y} \nonumber \\
    & \qquad + \left|\prob{\zt=m-y} - \prob{Z=m-y}\right|\prob{Y=y} \nonumber \\
    &\le\  \frac{2}{p^2}\left(1 - p^2\right)^{N/4}\left(\prob{\zt=m-y} + \prob{Y=y}\right),
\end{align}
and therefore
\begin{align} 
    \bigg|\prob{\yt=y}&\prob{\zt=m-y} - \prob{Y=y}\prob{Z=m-y}\bigg| \nonumber \\
    &\le\ \frac{2}{p^2}\left(1 - p^2\right)^{N/4}\left(\sum_{y=0}^{m}{\prob{\zt=m-y}} + \sum_{y=0}^{m}{\prob{Y = y}}\right) \nonumber \\
    &\le\ \frac{4}{p^2}\left(1 - p^2\right)^{N/4}.\label{eq:prob sums almost equal}
\end{align}

Now, \eqref{eq:wt eq m prob}, \eqref{eq:prob Wm prob Wtm almost equal}, and \eqref{eq:prob sums almost equal}, together with the triangle inequality, imply the desired result.
\end{proof}

\section{Exponential Bounds on Missing Many Summands}\label{sec: exponential bounds on missing summands}
This section presents analogs of \cite{LMO}'s equations (4.6) and (4.14) for arbitrary $p$. Our derivations start with different intuitions and follow a different approach, but interestingly, we arrive at largely the same bounds.

\subsection{A bound on the $k$\textsuperscript{th} moment of $Y$}
From \eqref{eq:prob n small in A + A}, the probability of non-inclusion of $n$ into $A + A$ is at most $(1 - p^2)^{n+1}$. The probability of multiple numbers being non-included is less than or equal to the probability of each individual non-inclusion. Therefore, for integers $0 \le n_1, \dots, n_k \le N$,
\begin{equation}\prob{n_1, \dots, n_k \notin A + A} \ \le\ \left(\sqrt{1 - p^2}\right)^{1 + \max\set{n_1, \dots, n_k}}.\end{equation}

The number of tuples $(n_1, \dots, n_k)$ with $\max\{n_1, \dots, n_k\} \le n$ is $(n + 1)^k$. The expectation value of $Y^k$ can then be bounded as
\begin{align} \begin{split} \label{eq:bound on k-th moment}
    \e{Y^k} &\ \le\ \sum_{n_1=0}^{N}{\dots \sum_{n_k=0}^{N}{\left(\sqrt{1 - p^2}\right)^{1 + \max\set{n_1, \dots, n_k}}}}\\
            &\ <\ \sum_{n=0}^{\infty}{\left((n + 1)^k - n^k\right)\left(\sqrt{1 - p^2}\right)^{n+1}}\\
            &\ \le\ k \sum_{n=1}^{\infty}{n^{k-1} \left(\sqrt{1 - p^2}\right)^n} \ <\ \frac{2 k!}{\alpha^{k}},
\end{split} \end{align}
where
\begin{equation}\label{eq:alpha defn}\alpha \ \coloneqq\ \log{\frac{1}{\sqrt{1 - p^2}}} \ =\ \left|\log{\sqrt{1 - p^2}}\right|.\end{equation}

The derivation assumed $k \ge 1$, but in fact \eqref{eq:bound on k-th moment} is also valid for $k = 0$ because $\e{Y^0} = 1$. Equation \eqref{eq:bound on k-th moment} is a rather crude bound, but it is tight enough for the application of Chernoff's inequality. For $p = 1/2$, $\sqrt{1 - p^2} \approx 0.87$ and $\alpha \approx 0.14$. Corollary \ref{corr:better bound on k-th moment} gives a slightly better bound on the $k$\textsuperscript{th} moment, one with an \dq{$\alpha$} of $-\log{(\varphi/2)} \approx 0.21$.

\subsection{Applying the Chernoff Bound}
Whenever $|t| < \alpha$, the moment generating function $M(t)$ is bounded by
\begin{align} \begin{split} \label{eq:bound on moment generating function}
    M(t) &\ =\ \e{e^{t Y}} \ =\ \sum_{k=0}^{\infty}{\frac{\e{Y^k} t^k}{k!}} \ \le\ \sum_{k=0}^{\infty}{\frac{\frac{2 k!}{\alpha^k} t^k}{k!}}\\
    &\ =\ 2 \sum_{k=0}^{\infty}{\left(\frac{t}{\alpha}\right)^k} \ =\ \frac{2}{1 - t/\alpha}.
\end{split} \end{align}

By the Chernoff bound,
\begin{equation}\prob{Y \ge n} \ \le\ \inf_{t > 0}{\{M(t) e^{-t n}\}} \ =\ \inf_{t > 0}{\left\{\frac{2 e^{-t n}}{1 - t/\alpha}\right\}}.\end{equation}

For $n > 1 / \alpha$, the infimum $e^{-\alpha n + 1} n$ is attained at $t = \alpha - 1/n$. For $n \le 1 / \alpha$, the infimum $1$ is at $t = 0$ and the bound is trivial. Thus, for $n > 1 / \alpha$,
\begin{equation} \label{eq:chernoff bound}
    \prob{Y \ge n} \ \le\ 2 \, \alpha n \, e^{-\alpha n + 1} \ =\ O(n \, e^{-\alpha n}).
\end{equation}

Recall that $\alpha = -\log{\sqrt{1 - p^2}}$.

\probabilityExponentialBoundCorollary

For $p = 1/2$, $\sqrt{1 - p^2} \approx 0.87$, $\alpha \approx 0.14$, and $1/\alpha \approx 6.95$, and the bound starts being valid at $n = 7$. Corollary \ref{corr:better bound on prob of missing many} gives a slightly better bound, $O\left((0.81 + \eps)^n\right)$, which corresponds to an \dq{$\alpha$} of $0.21$. We note $\prob{Y \ge m}$ cannot be bounded tighter than exponential. If $0, \dots, n/2$ are missing from $A$, then $0, \dots, n$ are missing from $A + A$. Therefore, for even $n$,
\begin{equation} \label{eq:missing more than n lower bound}
    \prob{Y \ge n} \ \ge\ \prob{0, \dots, n \notin A + A} \ \ge\ \prob{0, \dots, n/2 \notin A} \ =\ (1 - p)^{n/2}.
\end{equation}

\begin{figure}[h] 
\centering
\includegraphics[width=0.78\textwidth]{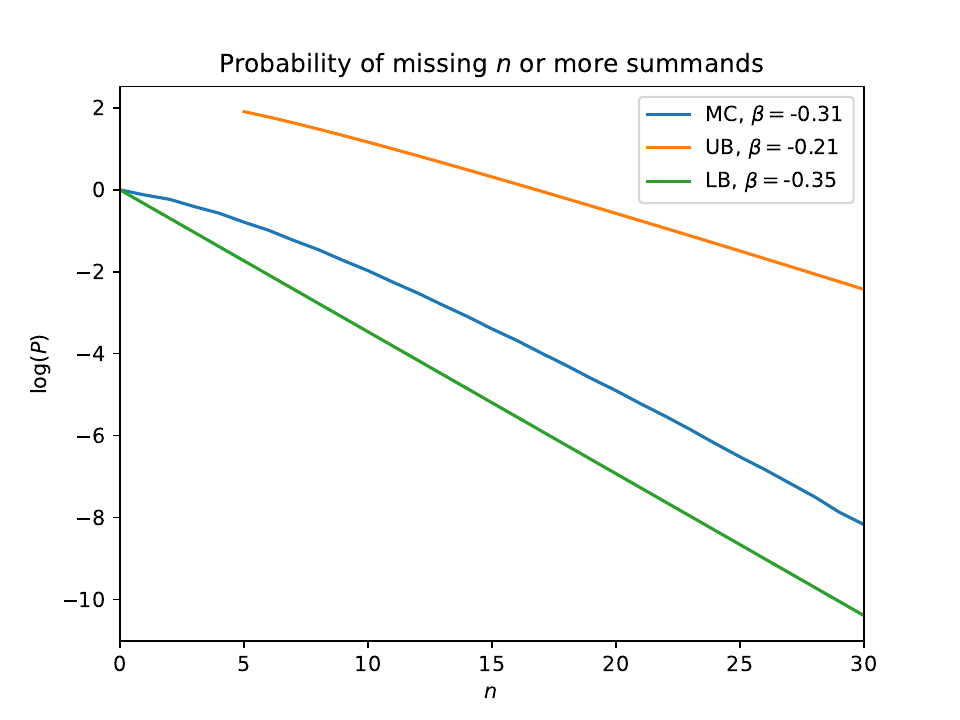}
\caption{Probabilities of missing more than $n$ summands on the left fringe for $p = 1/2$ and $N = 200$. Monte Carlo simulation (MC) with $1 \times 10^6$ trial runs, theoretical upper bound (UB) from Corollary \ref{corr:better bound on prob of missing many}, and theoretical lower bound (LB) from \eqref{eq:missing more than n lower bound} are shown. Here, $\beta$ is the slope of $\log(\prob{Y \ge n})$ as $n$ increases, equal to $\log{\sqrt{1 - p}}$ for LB and $\log{\lambda_1}$ for UB (see eq. \ref{eq:defn lambda12}). For MC, $\beta$ is estimated numerically by fitting a least-squares line, which entails not only random error but also systematic error because the original curve seems to be concave.}\label{fig:missing more than n}
\end{figure}

Corollary \ref{corr:better bound on prob of missing many} and \eqref{eq:missing more than n lower bound} together establish an approximate decay rate for $\prob{Y \ge n}$. Bounded above and below by two exponential functions, $\prob{Y \ge n}$ must itself be \dq{approximately exponential.} Figure \ref{fig:missing more than n} gives the result of Monte Carlo simulation ($1 \times 10^6$ trial runs) with $p = 1/2$ and $N = 200$, together with theoretical upper and lower bounds given by \eqref{eq:chernoff bound} and Corollary \ref{corr:better bound on prob of missing many}, respectively. The outcome of simulations is close to the lower bound and somewhat far from the upper bound.

By Theorem \ref{thm:W convolution}, \eqref{eq:chernoff bound} and \eqref{eq:missing more than n lower bound} about the probabilities of missing many summands from the left-fringe range $0, \dots, N$ translate into similar statements about missing many summands from the range $0, \dots, 2N$ of the entire sumset $A + A \subseteq \{0, \dots, 2N\}$.

\section{Probability of Missing Two Summands}\label{sec:twosummands}
This section provides a simplified expression for Proposition 3.5 from \cite{CKLMMSSX}, an exact formula for the probability of non-inclusion of two given numbers into $A + A$. The derivation of this simplified formula uses similar ideas as in the graph theoretic framework used in \cite{CKLMMSSX} and \cite{LMO}, but modified slightly. As such, we have opted to relegate most proofs to Appendix \ref{sec:appendix}. 

The probability $\prob{m, n \notin A + A}$ exhibits exponential decay in both $m$ and $n$ at a rate of about $\varphi/2 \approx 0.81$ for $p = 1/2$ (Corollary \ref{corr:missing two upper bound}), although the exact rate of decay depends on the ratio $l = \ceil{\frac{n + 1}{m - n}}$. The probability also depends on the parities of $m, n, l$. The numbers $a_k$ are introduced in Definition \ref{defn:a_k}, and is the same $a_k$ as in Lemma 3.1 of \cite{CKLMMSSX}. They start with $a_1 = 1$ and decay exponentially.

\begin{restatable}{lemma}{aKRecurrence} \label{lemma:a_k recurrence}
We have $a_1 = 1$, $a_2 = 1 - p^2$, and for $3 \le k \le N$, $a_k$ is given by the recurrence
\begin{equation}a_k \ =\ (1 - p)a_{k-1} + p(1 - p)a_{k-2}.\end{equation}
\end{restatable}

We now reformulate Prop 3.5 of \cite{CKLMMSSX}.

\begin{restatable}{proposition}{twoNoninclusionProbability} \label{prop:prob m, n notin sumset}
The probability of non-inclusion of two nonnegative integers $n < m \le N$ into $A + A$ is
$$\prob{m, n \notin A + A} \ =\ a_{2l+2}^{\floor{d_1/2}} a_{2l}^{\floor{d_2/2}} \begin{cases}
    1 & \text{$s = (1, 1, 0)$ or $(1, 1, 1)$}\\
    (1 - p) a_l & \text{$s = (1, 0, 1)$ or $(0, 1, 0)$}\\
    (1 - p) a_{l-1} & \text{$s = (1, 0, 0)$ or $(0, 1, 1)$}\\
    (1 - p)^2 a_l a_{l-1} & \text{$s = (0, 0, 0)$ or $(0, 0, 1)$},
\end{cases}$$
where $l = \ceil{\frac{n + 1}{m - n}}$ is the \dq{degree of twistedness} as in \eqref{eq:defn l}, $d_1$ and $d_2$ count the number of integers in $n+1, \dots, m$ greater than or equal to and less than the threshold $l(m - n)$, respectively,
\begin{equation}d_1 \ =\ (m+1) - l(m-n), \quad d_2 \ =\ l(m - n) - (n + 1),\end{equation}
and $s \coloneqq (m, n, l) \mod 2$ encodes the parities of $m, n, l$.
\end{restatable}

Define
\begin{equation} \label{eq:defn lambda12}
    \lambda_1 \ \coloneqq\ \frac{1 - p + \sqrt{(1 - p)(1 + 3p)}}{2}, \;\;\; \lambda_2 \ \coloneqq\ \frac{1 - p - \sqrt{(1 - p)(1 + 3p)}}{2}.
\end{equation}
One can check that $a_k$ is given by
\begin{equation} \label{eq:a_k binet}
    a_k \ =\ \frac{1 - \lambda_2}{\lambda_1 - \lambda_2} \lambda_1^k + \frac{-(1 - \lambda_1)}{\lambda_1 - \lambda_2} \lambda_2^k.
\end{equation}
In particular, since $\lambda_2 < \lambda_1$,
\begin{equation} \label{eq:a_k less than}
    \lambda_1^k \ <\ a_k \ \le\ \lambda_1^{k-1}.
\end{equation}

Proposition \ref{prop:prob m, n notin sumset} now gives us the following lemma.

\begin{restatable}{lemma}{missingTwoUpperBound}
\label{corr:missing two upper bound}
For $l \ge 2$, $(1 - p) a_{l - 1} \le \sqrt{a_{2l}}$. Moreover, since $a_k \le \lambda_1^{k-1}$, we have $\prob{m, n \notin A + A} \le a_{2l + 2}^{d_1/2} a_{2l}^{d_2/2}$ and $\prob{m, n \notin A + A} \le \lambda_1^{1 + \frac{m + n}{2}}$.
\end{restatable}

\begin{corollary} \label{corr:better bound on k-th moment}
The $k$\textsuperscript{th} moment of the number of missing summands on the left fringe is bounded by
\begin{equation}\e{Y^k} \ \le\ \e{Y} + \frac{2 k!}{\lambda_1 |\log{\lambda_1}|} \ \le\ \left(\frac{2}{p^2} - \frac{1}{p} - 1\right) + \frac{2 k!}{\lambda_1 |\log{\lambda_1}|}.\end{equation}
\end{corollary}

\betterProbabilityBoundCorollary

\begin{proof}
By Corollary \ref{corr:better bound on k-th moment}, the moment generating function is bounded by
\begin{equation}M(t) \ \le\ \left(\frac{2}{p^2} - \frac{1}{p} - 1\right) e^t + \frac{2/\lambda_1}{1 - t / \alpha'}.\end{equation}

By Chernoff's inequality, $\prob{Y \ge n}$ is less than or equal to $M(t) e^{-n t}$ for every $t > 0$ where $M(t)$ is defined. In particular, we may take $t = \alpha' - 1/n$.
\end{proof}

\section{The Second Moment}\label{sec: second moment}
There is little hope that $\prob{m, n \notin A + A}$, as given in Proposition \ref{prop:prob m, n notin sumset}, can be summed to a closed-form expression. To simplify, we take the limit $N \to \infty$. For small $p$, many summands will be missing from $A + A$, so we expect $\e{Y^2}$ to blow up to infinity as $p$ gets small. Since $\e{Y} = 2/p^2 - 1/p - 1$ in the limit $N \to \infty$, we hypothesize that $\e{Y^2}$ can be likewise expressed as a Laurent series in $p$. We find an exact (but not closed-form) expression for $\lim_{N\to\infty}{\e{Y^2}}$ and show that the leading term in its asymptotic expansion near $p = 0$ is $4 p^{-4}$.

Define
\begin{equation} \label{eq:defn S}
    S \ \coloneqq\ \lim_{N \to \infty}{\sum_{m=0}^{N}{\sum_{n=0}^{m-1}{\prob{m, n \notin A + A}}}},
\end{equation}
where $\prob{m, n \notin A + A}$ is given by Proposition \ref{prop:prob m, n notin sumset}, so that
\begin{align} \begin{split} \label{eq:eY2 in terms of S}
    \lim_{N \to \infty}{\e{Y^2}} \ =\ \lim_{N \to \infty}{\e{Y}} + 2 S \ =\ \left(\frac{2}{p^2} - \frac{1}{p} - 1\right) + 2 S.
\end{split} \end{align}

The difficulty in summing $\prob{m, n \notin A + A}$ is that it depends on $l$, the degree of twistedness defined in \eqref{eq:defn l}, but other than that, $\prob{m, n \notin A + A}$ is pretty much a geometric series. Define
\begin{equation} \label{eq:defn UV}
    U_l \ \coloneqq\ a_{2l}^l / a_{2l+2}^{l-1}, \quad V_l \ \coloneqq\ a_{2l+2}^l / a_{2l}^{l+1},
\end{equation}
so that for $n < m \le N$,
$$\prob{m, n \notin A + A} \ =\ U_l^{\frac{m+1}{2}} V_l^{\frac{n+1}{2}} \begin{cases}
    1 & \text{$s = (1, 1, 0)$ or $(1, 1, 1)$}\\
    (1 - p) a_l & \text{$s = (1, 0, 1)$ or $(0, 1, 0)$}\\
    (1 - p) a_{l-1} & \text{$s = (1, 0, 0)$ or $(0, 1, 1)$}\\
    (1 - p)^2 a_l a_{l-1} & \text{$s = (0, 0, 0)$ or $(0, 0, 1)$}.
\end{cases}$$

We break up the plane $\Z_{\ge 0}^2$ into \dq{wedges} of fixed $l$. For nonnegative integers $m$ and $l'$, define
\begin{equation}
    g_{m,l'} \ \coloneqq\ \left\lfloor\frac{(m + 1)(l'-1)}{l'}\right\rfloor,
\end{equation}
so that $l \ge l'$ (recall that $l \coloneqq \ceil{\frac{n+1}{m-n}}$) if and only if $n \ge g_{m,l'}$. For each $m$ and $l'$, the collection of $n$'s such that $l$ is equal to $l'$ is precisely $g_{m,l'}, \dots, g_{m,l'+1}-1$. Therefore,
\begin{equation}
    S \ =\ \lim_{N\to\infty}{\sum_{l'=1}^{\infty}{\sum_{m=0}^{\infty}{\sum_{n=g_{m,l'}}^{g_{m,l'+1}-1}{\prob{m, n \notin A + A}}}}}.
\end{equation}

We may define
\begin{align} \begin{split} \label{eq:defn S^ij_l}
    S^{11}_l &\ \coloneqq\ \sum_{\text{$m$ odd}}{\sum_{\text{$n$ odd}}{U_l^{\frac{m+1}{2}} V_l^{\frac{n+1}{2}}}},\\
    S^{00}_l &\ \coloneqq\ \left((1 - p)^2 a_l a_{l-1}\right) \sum_{\text{$m$ even}}{\sum_{\text{$n$ even}}{U_l^{\frac{m+1}{2}} V_l^{\frac{n+1}{2}}}},\\
    S^{10}_l &\ \coloneqq\ \left(\begin{cases}
        (1 - p) a_l & \text{$l$ odd}\\
        (1 - p) a_{l-1} & \text{$l$ even}
    \end{cases}\right) \sum_{\text{$m$ odd}}{\sum_{\text{$n$ even}}{U_l^{\frac{m+1}{2}} V_l^{\frac{n+1}{2}}}}\\
    S^{01}_l & \ \coloneqq\ \left(\begin{cases}
        (1 - p) a_{l-1} & \text{$l$ odd}\\
        (1 - p) a_l & \text{$l$ even}
    \end{cases}\right) \sum_{\text{$m$ even}}{\sum_{\text{$n$ odd}}{U_l^{\frac{m+1}{2}} V_l^{\frac{n+1}{2}}}},
\end{split} \end{align}
so that
\begin{equation}S \ =\ \sum_{l=1}^{\infty}{\left(S^{11}_l + S^{00}_l + S^{10}_l + S^{01}_l\right)}.\end{equation}

(A small abuse of notation: instead of summing over a dummy variable we sum over $l$) Each one of the sums in \eqref{eq:defn S^ij_l} is almost of the form $\sum_{i=0}^{\infty}{\sum_{j=0}^{i-1}{\alpha^i \beta^j}}$. For example,
\begin{align} \begin{split}
    S^{11}_l \ =\ \sum_{\text{$m$ odd}}{\sum_{\text{$n$ odd}}{U_l^{\frac{m+1}{2}} V_l^{\frac{n+1}{2}}}} &\ =\ \sum_{b=0}^{\infty}{U_l^{b+1} \sum_{a=g_{bl}}^{g_{b,l+1}-1}{V_l^{a+1}}}\\
        &\ =\ \frac{U_l V_l}{1 - V_l} \left(\sum_{b=0}^{\infty}{U_l^b V^{g_{bl}}} - \sum_{b=0}^{\infty}{U_l^b V_l^{g_{b,l+1}}}\right),
\end{split} \end{align}
where we used the substitution $m = 2 b + 1$ and $n = 2 a + 1$, and also the fact that $\ceil{\frac{(2a+1)+1}{(2b+1)-(2a+1)}} = \ceil{\frac{a+1}{b-a}}$. The exponent $g_{b,l+1}$ grows linearly with $b$, up to occasional corrections having to do with the floor function. Lemma \ref{lemma:floor-geometric sum}, which is proven below, allows us to evaluate such \dq{floor-geometric sums}. The result is that
\begin{equation}S^{11}_l \ =\ \sum_{l=1}^{\infty}{\frac{a_{2l}}{(1 - a_{2l+2})(1 - a_{2l})}}.\end{equation}

The other sums can be likewise evaluated, leading to
\begin{equation}
    S \ =\ \sum_{l=1}^{\infty}{\frac{a_{2l} + (1-p)a_{l-1} + (1-p)a_l a_{2l} + (1-p)^2 a_l a_{l-1}}{(1 - a_{2l+2})(1 - a_{2l})}} - \left(\frac{2}{p^2} - \frac{1}{p} - 1\right).
\end{equation}

From \eqref{eq:eY2 in terms of S}, we deduce the following.

\begin{proposition}
We have
\begin{align} \begin{split} \label{eq:second moment}
    \lim_{N\to\infty}{\e{Y^2}} \ = \ & - \left(\frac{2}{p^2} - \frac{1}{p} - 1\right) +\\
        &+2 \sum_{l=1}^{\infty}{\frac{a_{2l} + (1-p) a_{l-1} + (1-p)a_l a_{2l} + (1-p)^2 a_l a_{l-1}}{(1 - a_{2l+2})(1 - a_{2l})}}.
\end{split} \end{align}
\end{proposition}

This is an exact expectation value of the square of the number of missing summands in the left fringe as $N \to \infty$. The summands in \eqref{eq:second moment} decay exponentially, since $a_k \le \lambda_1^{k-1}$. Equation \eqref{eq:second moment} reduces to \cite{LMO}'s Theorem 1.5 for $p = 1/2$ and also works for general $p$. We achieved significant simplification over \cite{LMO} Theorem 1.5, in the sense that breaking the plane into \dq{wedges} and summing over these wedges instead of coordinate pairs converted a double infinite sum into a single infinite sum and removed floor functions and parity dependence.

\begin{figure}[h] 
\centering
\includegraphics[width=0.95\textwidth]{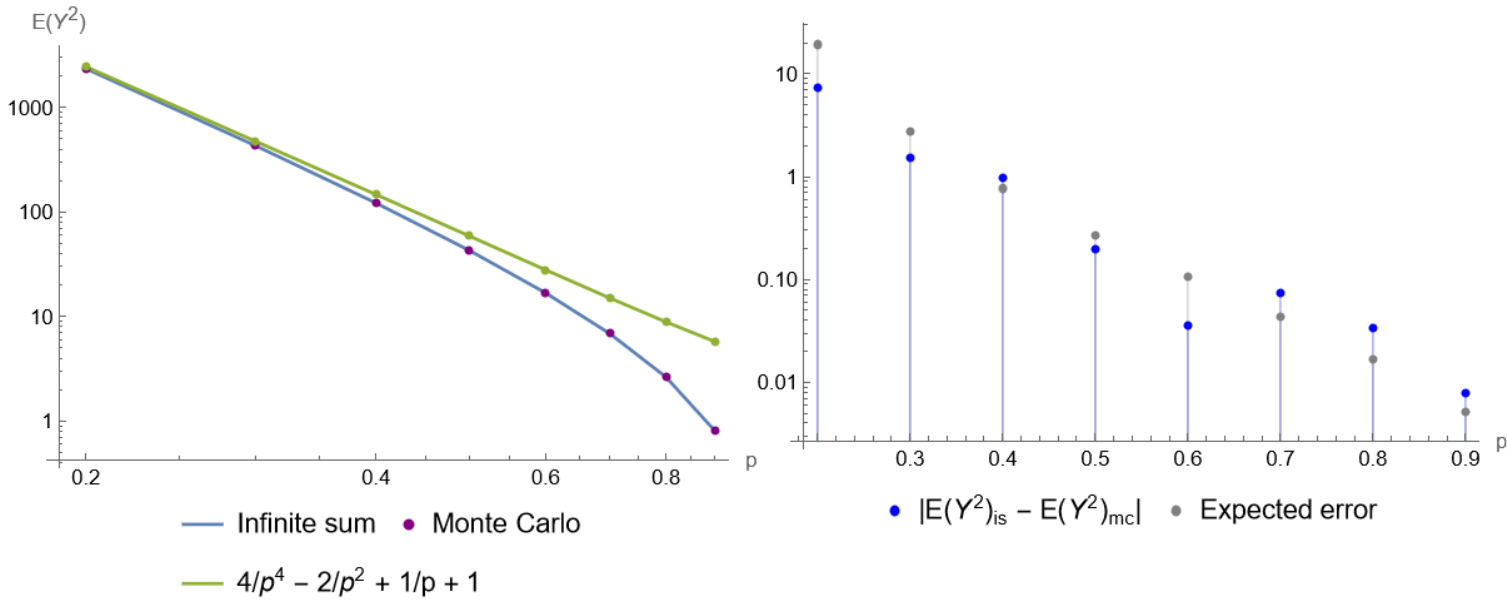}
\caption{Second moment $\e{Y^2}$ of the number of missing summands on the left fringe: theoretical prediction  \eqref{eq:second moment}, Monte Carlo values, and the $4 p^{-4}$ approximation \eqref{eq:2p-4 approx} (left); discrepancies between Monte Carlo values $\e{Y^2}_{\rm mc}$ and the infinite sum $\e{Y^2}_{\rm is}$ from \eqref{eq:second moment}, together with what we expect the discrepancies to be based on simulation size and $N$ (right).}\label{fig:second moment numerics}
\end{figure}

Figure \ref{fig:second moment numerics} compares the result of numerically evaluating \eqref{eq:second moment} and of Monte Carlo simulations with $M = 1 \times 10^5$ runs for some values of $p$ and for $N = 400$, as well as the $4 p^{-4}$ approximation which will be discussed in the next section. Expected errors were calculated as the square root of random error squared plus systematic error squared, where random error due to the finite size of the Monte Carlo simulation was simplistically taken to be $\Delta \e{Y^2} = 2 \e{Y^2} / \sqrt{M}$, and systematic error due to the finite $N$ used in the simulation was estimated as
\begin{align} \begin{split} \label{eq:second moment limit remainder estimate}
    \e{Y^2}_{\infty} - \e{Y^2}_N &\ =\ \sum_{\max\set{m,n} \ge N+1}{\prob{m, n \notin A + A}}\\
                                &\ \le\ \sum_{n=N+1}^{\infty}{((n+1)^2 - n^2) \left(\sqrt{1 - p^2}\right)^{n+1}}.
\end{split} \end{align}

The log-log plot is almost linear with a slope of $-4$, and the approximation $\e{Y^2} \approx 4/p^4 - 2/p^2 + 1/p + 1$ \eqref{eq:2p-4 approx} appears very accurate at small $p$. The fact that discrepancies between Monte Carlo and \eqref{eq:second moment} shrink together with expected errors serves as independent evidence that the complicated algebra leading up to \eqref{eq:second moment} was probably correct.

\begin{lemma} \label{lemma:floor-geometric sum}
For numbers $\alpha, \beta$ with $|\alpha|, |\beta| < 1$ and integers $0 \le k < l$,
\begin{equation}\sum_{n=0}^{\infty}{\alpha^n \beta^{\floor{\frac{(l-1)n + k}{l}}}} \ =\ \frac{1}{1 - \alpha \beta}\left(1 + \frac{\alpha^{k+1}\beta^k(1 - \beta)}{1 - \alpha^l \beta^{l-1}}\right).\end{equation}
\end{lemma}
\begin{proof}
As $n$ increases, $(l-1)n + k$ modulo $l$ goes as $k, k-1, \dots, 0$, then jumps to $l-1$ and settles into a cyclic mode $l-1, l-2, \dots, 0, l-1, l-2, \dots, 0$. As $n$ goes to $n+1$, the quantity $\floor{\frac{(l-1)n + k}{l}}$ increases by $1$ if $(l-1)n + k$ modulo $l$ is nonzero and stays the same if $(l-1)n + k$ modulo $l$ is zero. Motivated by this, we make the observation that for $0 \le n \le k$,
\begin{equation}\floor{\frac{(l-1)n + k}{l}} \ =\ n + \floor{\frac{k-n}{l}} = n,\end{equation}
and that if we write $n$ as $l s + d + k$ for some integers $s, d$ with $1 \le d \le l$, then
\begin{equation}\floor{\frac{(l-1)n + k}{l}} \ =\ \floor{(l-1)s+k + \frac{(l-1)d}{l}} \ =\ (l-1)s + k + d - 1.\end{equation}

If we split the sum in question into the regions $[0, k]$ and $[k+1, \infty)$ and apply the two equations above, the sum becomes a collection of geometric series each of which is easy to evaluate:
\begin{align} \begin{split}
\sum_{n=0}^{\infty}{\alpha^n \beta^{\floor{\frac{(l-1)n + k}{l}}}} &\ =\ \sum_{n=0}^{k}{\alpha^n \beta^n} + \sum_{s=0}^{\infty}{\sum_{d=1}^{l}{\alpha^{ls + d + k} \beta^{(l-1)s + d + k - 1}}}\\
    &\ =\ \sum_{n=0}^{k}{\left(\alpha \beta\right)^n} + \alpha^k \beta^{k-1} \left(\sum_{s=0}^{\infty}{\left(\alpha^l \beta^{l-1}\right)^s}\right)\left(\sum_{d=1}^{l}{\left(\alpha \beta\right)^d}\right).
\end{split} \end{align}
\end{proof}

\section{Concentration of $Y$ and the Asymptotics of the Second Moment}\label{sec: concentration and asymptotics}
The first moment $\e{Y}$ in the limit of large $N$ is $2p^{-2}$ to leading order in $p$. Therefore, we might expect the second moment $\e{Y^2}$ in this limit to grow as $(2p^{-2})^2 = 4p^{-4}$. We prove that that is indeed the case. 

Moreover, we show that the variance $\text{Var}{(Y)}$ is asymptotically strictly less than $p^{-4}$ and that therefore for small $p$, the number of missing summands $Y$ in the left fringe is concentrated around the mean $2/p^2 - 1/p - 1$. By Theorem \ref{thm:W convolution}, this translates into saying that the total number of missing summands is also concentrated.

\subsection{$\e{Y^2}$ to Leading Order}
Recall that \eqref{eq:second moment} provides an exact formula for $\lim_{N\to\infty}{\e{Y^2}}$ in terms of an infinite series whose value is hard to compute. Since we expect this infinite series to behave as $p^{-4}$, for $0 < p < 1$ and $L$ a positive integer, let us define $f(p, L)$ to be $p^4$ times its $L$\textsuperscript{th} partial sum,
\begin{equation} \label{eq:defn fpL}
    f(p, L) \ \coloneqq\ \sum_{l=1}^{L}{\frac{p^4 \left(a_{2l} + (1-p)a_{l-1} + (1-p)a_l a_{2l} + (1-p)^2 a_l a_{l-1}\right)}{(1 - a_{2l+2})(1 - a_{2l})}},
\end{equation}
where the $a_k$'s are to be understood as depending on this new value of $p$ as in Lemma \ref{lemma:a_k recurrence} or \eqref{eq:a_k binet}, so that we might compute $\lim_{p \to 0}{\e{ Y^2} / p^{-4}}$ as twice the limit of $f(p, L)$ as $L$ goes to infinity and $p$ to zero. The limit $L \to \infty$ is to be taken first and $p \to 0$ second. The other order of limits, when $p$ is first taken to zero and then $L$ to infinity, is much easier to evaluate. We evaluate this second order of limits and prove that since $f(p, L)$ is sufficiently convergent, the limits may be exchanged.

\begin{proposition} \label{prop:fpL pointwise to p}
We have $f(p, L)$ converges pointwise in $p$ to $\sum_{l=1}^{L}{\frac{4}{4l^2 - 1}}$.
\end{proposition}
\begin{proof}
From \eqref{eq:defn lambda12}, one can compute that $\lambda_1 = 1 - p^2 + O(\lambda^3)$, $\lambda_2 = O(p)$, $C_1 = 1 + p^2 + O(\lambda^3)$, $C_2 = O(p)$. Since $a_{2l} = C_1 \lambda_1^{2l} + C_2 \lambda_2^{2l}$,
\begin{equation}\lim_{p \to 0}{\frac{1 - a_{2l}}{p^2}} \ =\ 2l - 1.\end{equation}

As $p$ goes to zero, each summand in \eqref{eq:defn fpL} converges to $\frac{4}{4l^2 - 1}$. Since there are finitely many summands, the desired result follows.
\end{proof}

\begin{lemma} \label{lemma:converigng to one lower bound}
Let $\mu$ be a real number between $0$ and $1$, and let $R$ be a positive real number. For $x$ a real number between $0$ and $R$,
\begin{equation}1 - \mu^x \ \ge\ \frac{1 - \mu^R}{R} x.\end{equation}
\end{lemma}
\begin{proof}
$1 - \mu^x$ is a function with a strictly negative second derivative whose plot intersects the line $\frac{1 - \mu^R}{R} x$ at the two points $x = 0$ and $x = R$.
\end{proof}

\begin{proposition} \label{prop:fpL Cauchy in L}
We have $f(p, L)$ is uniformly Cauchy in $L$.
\end{proposition}
\begin{proof}
For any tolerance $\eps > 0$, which we may assume to be less than $1$, let $\eps' \coloneqq \eps / 64$, and let
\begin{equation}K \ \coloneqq\ \ceil{\frac{9 |\log{\eps'}|^2}{\eps'}}.\end{equation}

For integers $a, b$ with $K \le a < b$, the difference $f(p, b) - f(p, a)$ may be bounded as
\begin{equation} \label{eq:fpa - fpb}
    f(p, b) - f(p, a)\ <\ \sum_{l=a+1}^{\infty}{\frac{4 p^4 \lambda_1^{l-2}}{(1 - \lambda_1^{2l+1})(1 - \lambda_1^{2l-1})}},
\end{equation}
where we have used \eqref{eq:a_k less than}, $a_k \le \lambda_1^{k-1}$. One can check that $\lambda_1 = 1 - p^2 + O(p^3)$. The denominator in \eqref{eq:fpa - fpb} is approximately $(2l+1)(2l-1) p^4$ for small $l$ and approximately constant for large $l$. Define
\begin{equation}r \ \coloneqq\ \ceil{\frac{|\log{\eps'}|}{|\log{\lambda_1|}}} + 2\end{equation}
to quantify the boundary between the \dq{small $l$} and \dq{large $l$} regions, so that $\lambda^{l-2} \le \eps'$. For $l \ge r$, the numerator is small and the denominator is $1$ up to an $\eps'$-sized correction. Since $\lambda_1 \approx 1 - p^2$, the value of $r$ depends on $p$ as $O(1/p^2)$, reflecting the fact that $a_k$ takes longer to converge to zero if $p$ is small. The large-$l$ portion of the sum in \eqref{eq:fpa - fpb} can be bounded as
\begin{align} \begin{split} \label{eq:large l sum bound}
\sum_{l=r}^{\infty}{\frac{4 p^4 \lambda_1^{l-2}}{(1 - \lambda_1^{2l+1})(1 - \lambda_1^{2l-1})}} &\ <\ \sum_{l=r}^{\infty}{\frac{4 p^4 \lambda_1^{l-2}}{(1 - \eps')^2}}\\
    &\ =\ \frac{8 p^2}{(1 - \eps')^2} \frac{p^2/2}{1 - \lambda_1} \lambda_1^{r-2}\\
    &\ <\ 32 \eps' \ =\ \frac{\eps}{2}.
\end{split} \end{align}

We upper-bound the sum in \eqref{eq:fpa - fpb} for all $0 < p < 1$ by considering the large-$p$ and small-$p$ regimes, where \dq{largeness} of $p$ is defined in reference to $a$ and $\eps'$. When $r \le a+1$, which is the large-$p$ regime, the entire sum bounding $f(p, b) - f(p, a)$ in \eqref{eq:fpa - fpb} is itself bounded by the geometric series from  \eqref{eq:large l sum bound}, so we're done. Let us consider the small $p$ regime when $r > a + 1$ (and also $a \ge K$ by assumption). From the definition of $r$, and from the fact that $|\log{\lambda_1}| \ge 1 - \lambda_1 \ge p^2/2$ we have the inequalities
\begin{align} \begin{split}
    & p^2/2 \ \le\ |\log{\lambda_1}| \ \le\ \frac{|\log{\eps'}|}{r - 3} \ \le\ \frac{|\log{\eps'}|}{K - 1},\\
    & 2 r - 1 \ \le\ \frac{2 |\log{\eps'}|}{|\log{\lambda_1}|} + 5 \ \le\ \frac{3 |\log{\eps'}|}{|\log{\lambda_1}|} \ \le\ \frac{6 |\log{\eps'}|}{p^2},
\end{split} \end{align}
where we also use that $K \ge 6$. By Lemma \ref{lemma:converigng to one lower bound} we have that for $0 \le x \le 2 r - 1$,
\begin{equation}\frac{1 - \lambda_1^x}{p^2} \ \ge\ \frac{1 - \lambda_1^{2r-1}}{p^2(2r - 1)} x \ >\ \frac{1 - \eps'}{6|\log{\eps'}|} x \ \ge\ \frac{x}{12 |\log{\eps'}|},\end{equation}
and therefore the small-$l$ region of the sum in  \eqref{eq:fpa - fpb} may be bounded as
\begin{align} \begin{split} \label{eq:small l sum bound}
    \sum_{l=a+1}^{r-1}{\frac{4 p^4 \lambda_1^{l-2}}{(1 - \lambda_1^{2l+1})(1 - \lambda_1^{2l-1})}} &\ <\ \sum_{l=a+1}^{r-1}{\frac{4(12 |\log{\eps'}|)^2 \lambda_1^{l-2}}{(2l+1)(2l-1)}}\\
        &\ <\ 576 |\log{\eps'}|^2 \sum_{l=a+1}^{\infty}{\frac{1}{4l^2 - 1}}\\
        &\ =\ \frac{576 |\log{\eps'}|^2}{2 a + 1} \ <\ \frac{288 |\log{\eps'}|^2}{K} \ \le\ \frac{\eps}{2}.
\end{split} \end{align}

Combining  \eqref{eq:large l sum bound} and \eqref{eq:small l sum bound}, for the small-$p$ case we obtain the bound
\begin{align} \begin{split}
    \sum_{l=a+1}^{\infty}{\frac{4 p^4 \lambda_1^{l-2}}{(1 - \lambda_1^{2l+1})(1 - \lambda_1^{2l-1})}} &\ <\ \frac{\eps}{2} + \frac{\eps}{2} \ =\ \eps
\end{split} \end{align}
as desired. Since the threshold $K$ does not depend on $p$, $f(p, L)$ is uniformly Cauchy in $L$.
\end{proof}

Propositions \ref{prop:fpL pointwise to p} and \ref{prop:fpL Cauchy in L} allow us to do an exchange of limits:
\begin{multline}
    \lim_{p \to 0}{\sum_{l=1}^{\infty}{\frac{p^4 \left(a_{2l} + (1-p)a_{l-1} + (1-p)a_l a_{2l} + (1-p)^2 a_l a_{l-1}\right)}{(1 - a_{2l+2})(1 - a_{2l})}}} \ =\\
        =\ \lim_{p \to 0}{\lim_{L \to \infty}{f(p, L)}} \ =\ \lim_{L \to \infty}{\lim_{p \to 0}{f(p, L)}} \ =\ \sum_{l=1}^{\infty}{\frac{4}{4l^2 - 1}} \ =\ 2.
\end{multline}

Therefore, the infinite sum from \eqref{eq:second moment} is, to leading order, given by
\begin{equation} \label{eq:infinite sum leading order}
    \sum_{l=1}^{\infty}{\frac{a_{2l} + (1-p)a_{l-1} + (1-p)a_l a_{2l} + (1-p)^2 a_l a_{l-1}}{(1 - a_{2l+2})(1 - a_{2l})}} \ =\ \frac{2}{p^4} + \dots,
\end{equation}
where the \dq{$\cdots$} represents an error term $\delta(p)$ with $\lim_{p \to 0}{\delta(p) / p^{-4}} = 0$. So, the second moment $\e{Y^2}$ can be roughly approximated as
\begin{equation}
    \lim_{N \to \infty}{\e{Y^2}} \ \approx\ \frac{4}{p^4} - \frac{2}{p^2} + \frac{1}{p} + 1. \label{eq:2p-4 approx}
\end{equation}

This approximation was obtained by combining  \eqref{eq:second moment} and \eqref{eq:infinite sum leading order}. While it probably does not represent the correct asymptotic expansion of $\e{Y^2}$ to zeroth order, it seems like a reasonable approximation. Figure \ref{fig:second moment numerics} compares the results of  \eqref{eq:2p-4 approx} with those of  \eqref{eq:second moment} and of Monte Carlo simulations. What can be concluded  from  \eqref{eq:second moment} and \eqref{eq:infinite sum leading order} the following.

\secondMomentLeadingOrderProp

 By Theorem \ref{thm:W convolution}, for the second moment of the total number of missing summands in the large-$N$ limit we likewise have
\begin{equation} \label{eq:total second moment leading order}
    \lim_{N \to \infty}{\e{W^2}} \ =\ \lim_{N \to \infty}{\left(2 \e{Y^2} + 2\e{Y}^2\right)} \ =\ \frac{16}{p^4} + \cdots.
\end{equation}

The value of $N$ for which $\e{Y^2}$ and $\e{W^2}$ approach their $N \to \infty$ limits increases as $p$ becomes smaller. By \eqref{eq:second moment limit remainder estimate}, the value of $N$ needed for $\e{Y^2}$ to get to within a tolerance of $\eps$ of its $N \to \infty$ limit can be roughly approximated as
\begin{equation}N_{\eps} \ \approx\ \frac{8 |\log{p}|}{p^2} + \frac{2 |\log{(\eps/8)}|}{p^2}.\end{equation}

Equations \eqref{eq:second moment leading order} and \eqref{eq:total second moment leading order} concern the double limit $\lim_{p \to 0}{\lim_{N \to \infty}{p^4 (\dots)}}$, and the ordering of the limits matters. In fact, for any fixed $N$, $\e{Y^2}$ is just bounded by $(N + 1)^2$ and $\e{W^2}$ by $(2 N + 1)^2$ because the two tandom variables count missing summands in $0, \dots, N$ and $0, \dots, 2 N$, respectivey, so the limit $\lim_{N \to \infty}{\lim_{p \to 0}{p^4 (\dots)}}$ if the \dq{$\dots$} is substituted with $\e{Y^2}$ or $\e{W^2}$ is simply zero.

\subsection{Concentration of $Y$}
As $p$ goes to zero, the variance of $Y$ grows slower than the mean squared. That is, if we denote by $\e{Y}_\infty$ the limit of $\e{Y}$ as $N$ goes to infinity, and likewise for $\e{Y^2}$, they relate to each other in such a way that the variance is small,
\begin{equation} \label{eq:variance small}
    \lim_{p\to 0}{\frac{\e{Y^2}_\infty - \left(\e{Y}_\infty\right)^2}{\left(\e{Y}_\infty\right)^2}} \ =\ \lim_{p\to 0}{\left(\frac{p^4 \e{Y^2}_\infty}{4} - 1\right)} \ =\ 0,
\end{equation}
where we have used that $\e{Y}_\infty = 2 p^{-2} + o(p^{-2})$ and $\e{Y_2}_\infty = 4 p^{-4} + o(p^{-4})$,  \eqref{en:e y} and \eqref{eq:second moment leading order}. 

This means that for smaller and smaller $p$, if $N$ is taken sufficiently large (depending on $p$), the distribution of the number of missing summands on the left fringe $Y$ will be sharply concentrated around the mean $2/p^2 - 1/p - 1$. 

The Central Limit Theorem with Weak Dependence does not apply to $Y = \sum{X_n}$ because the $X_n$'s do not represent a stationary process, but \eqref{eq:variance small} does prove that for small $p$ and sufficiently large $N$, the random variable $Y / \e{Y}$ approaches a delta distribution. 

\begin{figure}[h]
\centering
\includegraphics[width=0.7\textwidth]{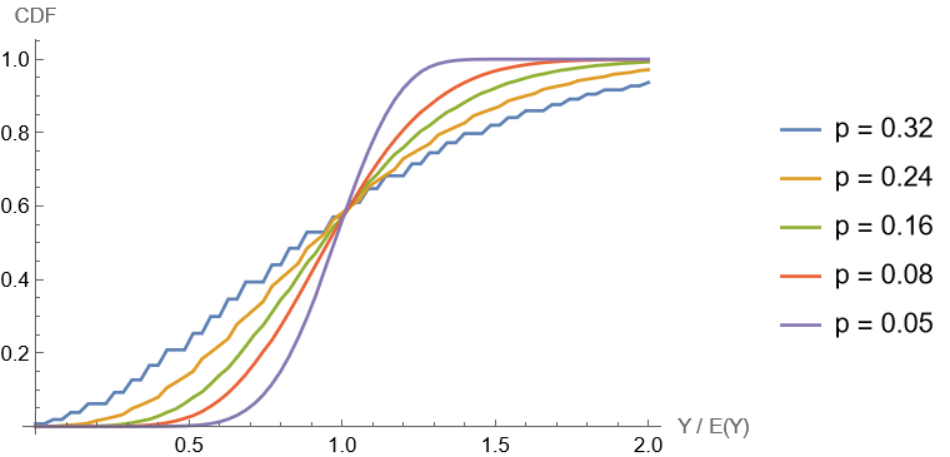}
\caption{The cumulative distribution function of $Y$, normalized by $\e{Y}$, for $N = 800$ and $p = 0.05, 0.08, 0.16, 0.24, 0.32$. (Monte Carlo simulation.)} \label{fig:cdf converging to step}
\end{figure}

Figure \ref{fig:cdf converging to step} demonstrates the cumulative probability distribution of $Y$ normalized by the mean $\e{Y}$. For $N = 800$ and each of several values of $p$, the following Monte Carlo simulation was run. A collection of $1 \times 10^6$ random sets $A$ were generated and the sumsets $A + A$ computed, and the mean number of missing summands across these subsets was obtained, and the number of sumsets missing less than a given fraction of the mean was counted. We note that the mean in each simulation was taken to be the sample mean and not $2/p^2 - 1/p - 1$, although the two are close. As $p$ gets smaller, the cumulative distribution function approaches a step function $\prob{Y \le y} \approx H(y / \e{Y} - 1)$.

\section{Future Work}\label{sec: future work}
We list a few natural directions for future work. 
\begin{enumerate}
    \item In \textsection \ref{sec: exponential bounds on missing summands}, the bounds for $\prob{n_1,\dots, n_k\notin A+A}$ and $\e{Y^k}$ are not optimal, as shown in Figure \ref{fig:missing more than n}. Does there exist a better upper bound for $\e{Y^k}$, and hence a better upper bound for $\prob{Y\geq n}$?
    \item In \textsection \ref{sec: second moment}, we take the limit $N\to \infty$ to obtain a closed-form expression for $\e{Y^2}$ in the limit. Does there exist a closed-form expression for $\e{Y^2}$ in general?
    \item In \textsection \ref{sec: concentration and asymptotics}, we isolate the leading term of $\e{Y^2}$ in the limit $N\to \infty$. Can lower order terms in this expansion be found? More generally, can an explicit expansion in terms of $1/p$ be found?
\end{enumerate}

These results only study the sumset $A+A$. The behavior of missing differences in the difference set $A-A$ in the limit as $N\to\infty$ has not been studied extensively. How can these results be modified to be applicable to difference sets?

\newpage

\appendix

\section{Appendix}\label{sec:appendix}
This section derives Proposition \ref{prop:prob m, n notin sumset}, an exact formula for the probability of non-inclusion of two given numbers into $A + A$, equivalent to but simpler than Proposition 3.5 from \cite{CKLMMSSX}. The probability $\prob{m, n \notin A + A}$ exhibits exponential decay in both $m$ and $n$ at a rate of about $\varphi/2 \approx 0.81$ for $p = 1/2$ (Corollary \ref{corr:missing two upper bound}), although the exact rate of decay depends on the ratio $l = \ceil{\frac{n + 1}{m - n}}$. The probability also depends on the parities of $m, n, l$. The numbers $a_k$ are introduced in Definition \ref{defn:a_k}, they start with $a_1 = 1$ and decay exponentially.

\twoNoninclusionProbability*

\begin{figure}[h] 
\caption{Snail diagram for $m = 17$, $n = 13$.}\label{fig:snail_diagram_loopless}
\centering
\includegraphics[width=0.9\textwidth]{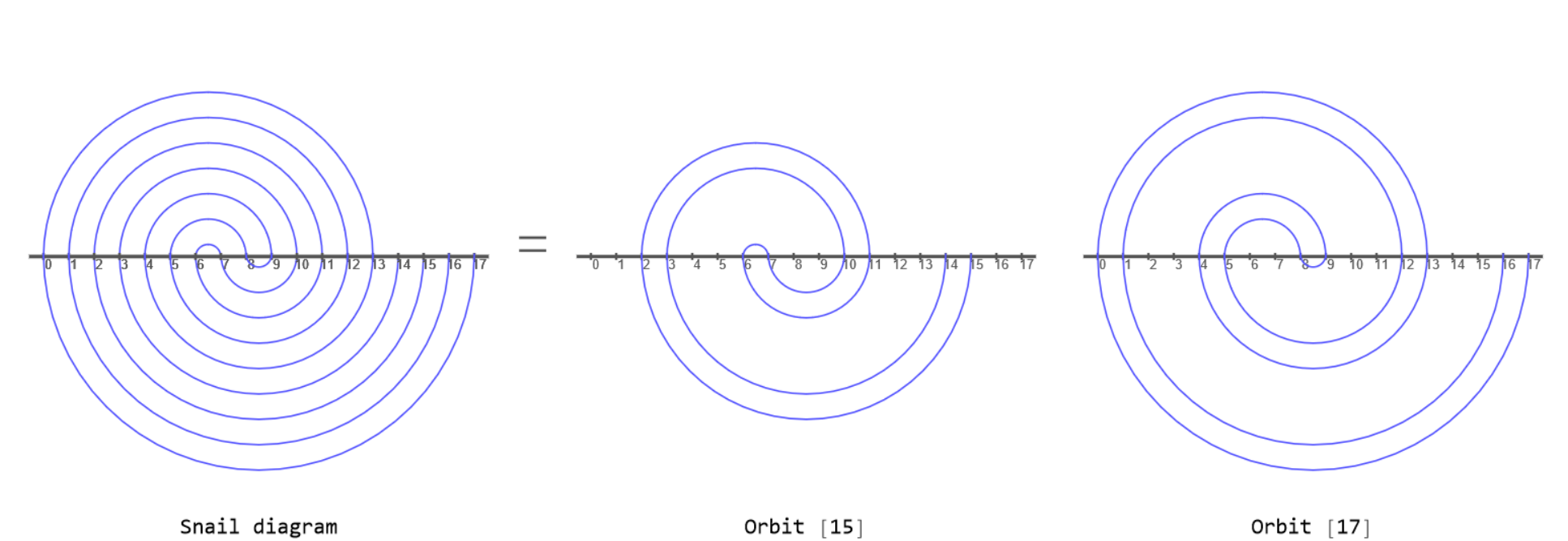}
\end{figure}

To visualize the condition for $m, n \notin A + A$, draw the number line. For every $0 \le x \le m$, connect $x$ with $m - x$ using a semicircle that goes below the number line, and for every $0 \le x \le n$, connect $x$ with $n - x$ using a semicircle that goes above the number line. Then, $m, n \notin A + A$ if and only if for every drawn semicircle, at least one of the endpoints lies in $A^c$. An example of such a diagram, which we call snail, is given in Figure \ref{fig:snail_diagram_loopless} for $m = 17$, $n = 13$. It is curious that the diagram splits in two disjoint spirals. In general, a diagram will split into a collection of spirals, and it is easy to calculate the lengths and counts of these spirals.

\begin{figure}[h] 
\caption{Snail diagram for $m = 17$, $n = 12$.}
\centering
\includegraphics[width=0.9\textwidth]{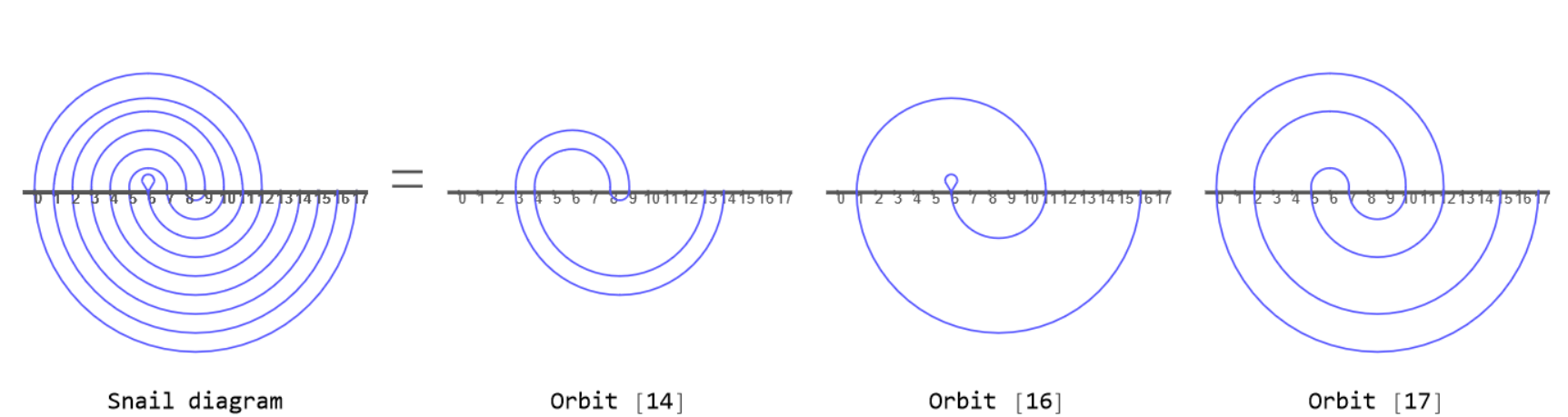}\label{fig:snail diagram loop}
\end{figure}

Spirals can either have two free ends, as in Figure \ref{fig:snail_diagram_loopless}, or have one free and one looped end, as in Figure \ref{fig:snail diagram loop}. We define orbits to formalize the pictorial notion of spirals. The orbits associated with Figure \ref{fig:snail_diagram_loopless} are $(15, 2, 11, 6, 7, 10, 3, 14)$ and $(17, 0, 13, 4, 9, 8, 5, 12, 1, 16)$, as well as the \dq{reverses} $(14, 3, 10, 7, 6, 11, 2, 15)$ and $(16, 1, 12, 5, 8, 9, 4, 13, 0, 17)$. The orbits associated with Figure \ref{fig:snail diagram loop} are $(14, 3, 9, 8, 4, 13)$, $(16, 1, 11, 6, 6)$, and $(17, 0, 12, 5, 7, 10, 2, 15)$, as well as the \dq{reverses} $(13, 4, 8, 9, 3, 14)$ and $(15, 2, 10, 7, 5, 12, 0, 17)$. The fact that $(16, 1, 11, 6, 6)$ ends with a repeated $6$ is called a \emph{loop}. Looped orbits don't have reverses, so the repeated number is always at the end and not at the start.

\begin{definition} \label{defn:orbit}
Let $n < r \le m$ be nonnegative integers. We wish to study the numbers one can reach from $r$ via the reflection maps $i \mapsto m - i$ and $i \mapsto n - i$ while staying in the nonnegative integers. Let $r_1 = r$, and for $t \ge 2$, let $r_t = m - r_{t-1}$ for $t$ even and $r_t = n - r_{t-1}$ for $t$ odd. Define the \emph{orbit} of $r$ as
\begin{equation}[r]_{mn} \ =\ (r_1, \dots, r_k),\end{equation}
where $k$ is the smallest positive integer such that $r_{k+1}$ is either negative or equals $r_{k-1}$. (Such a $k$ always exists because $r_{2p+1} = r - p(m - n)$ decreases with $p$.)
\end{definition}

One can check that
\begin{equation} \label{eq:orbit r}
r_t \ =\ \begin{cases}
    r - \frac{t-1}{2}(m - n) & \text{$t$ odd}\\
    m - r + \frac{t-2}{2}(m - n) & \text{$t$ even}
\end{cases}
\end{equation}
for all $t \ge 1$.

\begin{definition} \label{defn:chain}
A tuple $(x_1, \dots, x_k)$ of nonnegative integers is called a \emph{chain} of length $k$. This chain is said to be \emph{satisfied} if for every $1 \le i \le k-1$, either $x_i \notin A$ or $x_{i+1} \notin A$ (or both).
\end{definition}

\begin{lemma}[Collusion of chains] \label{lemma:chain collusion}
If chains $(r_1, \dots, r_k)$ and $(q_1, \dots, q_l)$ are such that $r_k = q_1$, then they are simultaneously satisfied if and only if the chain $(r_1, \dots, r_k, q_2, \dots, q_l)$ is satisfied.
\end{lemma}
\begin{proof}
The combined chain has the same adjacencies as the first two chains.
\end{proof}

\begin{theorem} \label{thm:two integers non-included orbits}
Two nonnegative integers $n < m$ are simultaneously non-included in $A + A$ if and only if every one of the chains $[r]_{mn}$ for $n < r \le m$ is satisfied.
\end{theorem}
\begin{proof}
If some orbit $[r]_{mn} = (r_1, \dots, r_k)$ is not satisfied, there exists some $i$ such that $r_i$ and $r_{i+1}$ are both in $A$. Then $r_i + r_{i+1}$ is in $A + A$, and by the definition of an orbit, $r_i + r_{i+1}$ is either $m$ or $n$.

Supose that $[r]_{mn}$ is satisfied for all $r$. Observe that $m$ is non-included in $A + A$ if and only if every one of the chains $(0, m), (1, m-1), \dots, (m, 0)$ is satisfied, and likewise for $n$. Define the set $C$ as
\begin{equation}C \ \coloneqq\ \set{(0, m), (1, m-1), \dots, (m, 0)} \sqcup \{(0, n), (1, n-1), \dots, (n, 0)\},\end{equation}
so that $m, n \notin A + A$ if and only if every one of the chains in $C$ is satisfied. We can consider $C$ modulo the equivalence relation $\sim$, where $(x, y) \sim (z, w)$ if and only if there exists a sequence of numbers $a_0, \dots, a_k$ with $a_0 \in \set{x, y}$, $a_k \in \set{z, w}$, and $(a_i, a_{i+1}) \in C$ for all $0 \le i < k$. Let $E$ be an equivalence class, and let
\begin{equation}
r \ =\ \max\set{x: \exists y \in \Z_{\ge 0}: \text{$(x, y) \in E$ or $(y, x) \in E$}}
\end{equation}
be the greatest number encountered in $E$. Then $r > n$, because if $r$ were less than or equal to $n$, $(r, n - r)$ would be in $E$, and $(m - (n - r), n - r)$ would be in $E$, but $r + (m - n) > r$.

Since $n - r$ is negative, $(r, n - r)$ or $(n - r, r)$ cannot be in $E$, and it has to be that $(r, m - r)$ and $(m - r, r)$ are in $E$. One can verify that if the orbit of $r$ is $[r]_{mn} = (r_1, \dots, r_k)$, then
\begin{equation}\set{(r_1, r_2), (r_2, r_3), \dots, (r_{k-1}, r_k)} \cup \set{(r_2, r_1), (r_3, r_2), \dots, (r_k, r_{k-1})}\end{equation}
is a set containing $(r, m - r)$ and closed under the discussed equivalence relation. Therefore, $E$ is equal to this set. By Lemma \ref{lemma:chain collusion}, all chains in $E$ are satisfied if and only if $[r]_{mn}$ is satisfied, which is true by assumption. Since $C$ is equal to the disjoint union of its equivalence classes, and every chain in every equivalence class is satisfied, $m$ and $n$ are non-included in $A + A$.
\end{proof}

\subsection{Loopless orbits}
\begin{lemma} \label{lemma:equiv m/2 n/2 mod m - n}
If $x$ and $n < m$ are integers, and if $x \equiv n/2 \pmod{m - n}$, then $m - x$ and $n - x$ are also equivalent to $n/2$;  if $x \equiv m/2 \pmod{m - n}$, then $m - x$ and $n - x$ are also equivalent to $m/2$.
\end{lemma}
\begin{proof}
For the first part, $m - n/2 = (m - n) + n/2 \equiv n/2$ and $n - n/2 = n/2$. For the second part, $m - m/2 = m/2$ and $n - m/2 = -(m - n) + m/2 \equiv m/2$.
\end{proof}

Throughout the rest of this subsection, we let $n < r \le m$ be nonnegative integers and $[r]_{mn} = (r_1, \dots, r_k)$ the orbit of $r$.

\begin{lemma} \label{lemma:loopless orbit}
If $r$ is not equivalent to $m/2$ or $n/2$ modulo $m - n$, then $r_1, \dots, r_k$ are all distinct. (We say that the orbit is \emph{loopless}.)
\end{lemma}
\begin{proof}
By  \eqref{eq:orbit r}, the odd entries of $r$ are strictly decreasing and the even entries are strictly increasing, so $r_s$ cannot equal to $r_t$ when $s \ne t$ have the same parity. To rule out the other case, suppose $t - s = 2 d + 1$ for some nonnegative integer $d$. If $s$ is odd and $t$ even, $r_s = m - r_{s+1}$ and $r_t = m - r_{t-1}$, and if $s$ is even and $t$ odd, $r_s = n - r_{s+1}$ and $r_t = n - r_{t-1}$; in either case, equality $r_s = r_t$ implies $r_{s+1} = r_{t-1}$. Applying this result inductively $d$ times, we find that $r_u = r_{u+1}$, where $u = s + d$. Therefore, either $r_u = n/2$ or $r_u = m/2$. By Lemma \ref{lemma:equiv m/2 n/2 mod m - n} and by induction, this would imply that $r$ is equivalent to either $m/2$ or $n/2$ modulo $m - n$, which is not the case.
\end{proof}

\begin{lemma} \label{lemma:loopless orbit length}
If $r$ is not equivalent to $m/2$ or $n/2$ modulo $m - n$, the length of its orbit is
\begin{equation}k \ =\ 2 \bigg\lceil \frac{r + 1}{m - n}\bigg\rceil.\end{equation}
\end{lemma}
\begin{proof}
By Lemma \ref{lemma:loopless orbit}, the $r_t$ are all distinct, so by Definition \ref{defn:orbit}, $k$ is simply the smallest positive integer such that $r_{k+1}$ is negative. By  \eqref{eq:orbit r}, the even entries $r_{2p+2}$ are all positive, so it has to be that $k + 1$ is odd. Thus,
\begin{equation}r_{k+1} \ =\ r - \frac{k}{2}(m - n) \ \le\ -1,\end{equation}
which is equivalent to saying that $\frac{k}{2} \ge \frac{r + 1}{m - n}$. By definition of the ceiling function, the smallest $k$ that satisfies this is $2 \ceil{\frac{r + 1}{m - n}}$.
\end{proof}
\begin{corollary} \label{corr:r_k k}
The endpoint $r_k$ satisfies $k \ =\ 2 \ceil{\frac{r_k + 1}{m - n}}$.
\end{corollary}
\begin{proof}
The orbit of $r_k$ is $[r_k]_{mn} = (r_k, \dots, r_1)$, the reverse of $[r]_{mn}$. It has the same length, so $k = 2 \ceil{\frac{r_k + 1}{m - n}}$.
\end{proof}

\subsection{Looped orbits}
\begin{lemma} \label{lemma:looped orbit length}
If $r$ is equivalent to $m/2$ or $n/2$ modulo $m - n$, the length of its orbit is
\begin{equation}k \ =\ \ceil{\frac{r + 1}{m - n}} + 1.\end{equation}
\end{lemma}
\begin{proof}
Suppose $r$ is equivalent to $m/2$. Write $r = m/2 + j(m - n)$ for some integer $j$. Note that $j$ is the unique integer satisfying $n < m/2 + j(m - n) \le m$. In particular, $j$ is the smallest integer satisfying $n + 1 \le m/2 + j(m - n)$, so we may write
\begin{equation}j \ =\ \ceil{\frac{n + 1 - m/2}{m - n}} \ =\ \ceil{\frac{m/2 + 1}{m - n}} - 1.\end{equation}

The orbit of $r$ includes $j$ repetitions of $x \mapsto m - x \mapsto x - (m - n)$ followed by one final $m$-reflection $m/2 \mapsto m - m/2$. Since there are $2 j + 1$ reflections done, the number of points $k$ is $2j + 2$. Then
\begin{equation}\ceil{\frac{r + 1}{m - n}} \ =\ j + \ceil{\frac{m/2 + 1}{m - n}} \ =\ 2j + 1 \ =\ k - 1,\end{equation}
which is what we wanted to show. Now suppose $r$ is equivalent to $n/2$. Write $r = n/2 + j(m - n)$ for some integer $j$. Note again that $j$ is the unique integer satisfying $n < n/2 + j(m - n) \le m$. In particular, $j$ is the smallest integer satisfying $n + 1 \le n/2 + j(m - n)$, so we may write
\begin{equation}j \ =\ \bigg\lceil \frac{n/2 + 1}{m - n}\bigg\rceil.\end{equation}

The orbit of $r$ includes $j$ repetitions of $x \mapsto m - x \mapsto x - (m - n)$. Since there are $2 j$ reflections done, the number of points $k$ is $2j + 1$. Then
\begin{equation} \label{eq:l l+1 inequality}
    \ceil{\frac{r + 1}{m - n}} \ =\ j + \ceil{\frac{n/2 + 1}{m - n}} \ =\ 2j \ =\ k - 1,
\end{equation}
which is what we wanted to show.
\end{proof}

\subsection{Orbit counts}
Throughout this subsection, $n < m$ are nonnegative integers and the \dq{degree of twistedness}
\begin{equation} \label{eq:defn l}
    l \ \coloneqq\ \ceil{\frac{n + 1}{m - n}}
\end{equation}
is a useful number that encodes how many twists the snail diagram for $m, n$ has. The calculational significance of $l$ is that for $n + 1 \le x \le l(m - n) - 1$,
\begin{equation}l \ \le\ \ceil{\frac{n + 2}{m - n}} \ \le\ \ceil{\frac{x + 1}{m - n}} \ \le\
\ceil{\frac{(l(m - n) - 1) + 1}{m - n}} \ =\ l,\end{equation}
and for $l(m - n) \le x \le m$,
\begin{equation}l + 1 \ =\ \ceil{\frac{l(m - n) + 1}{m - n}} \ \le\ \ceil{\frac{x + 1}{m - n}} \ \le\ \ceil{\frac{m + 1}{m - n}} \ =\ l + 1,\end{equation}
so that for $n + 1 \le x \le m$,
\begin{equation} \label{eq:ceil x + 1 over m - n}
    \ceil{\frac{x + 1}{m - n}} \ =\ \begin{cases}
        l & x < l(m - n)\\
        l + 1 & x \ge l(m - n).
    \end{cases}
\end{equation}

Let $a_1$ be the number of loopless orbits of length $2(l + 1)$ up to reversal (e.g., $(r_1, \dots, r_k)$ is declared equivalent to $(r_k, \dots, r_1)$), and $a_2$ the number of loopless orbits of length $2 l$ up to reversal. Then, the number of loopless orbits of length $2(l + 1)$ is $2 a_1$, and those of length $2 l$ is $2 a_2$. By Lemma \ref{lemma:loopless orbit length}, for the loopless orbits of length $2(l + 1)$, both endpoints have to be in $l(m - n), \dots, m + 1$, and for the loopless orbits of length $2 l$, both endpoints have to be in $n + 1, \dots, l(m - n) - 1$. Let $c_1$ be the number of integers in $l(m - n), \dots, m + 1$ that do not constitute the endpoint of a loopless orbit, and let $c_2$ be the number of integers in $n + 1, \dots, l(m - n) - 1$ that do not constitute the endpoint of a loopless orbit. Denote $d_1 \coloneqq (m + 1) - l(m - n)$ and $d_2 \coloneqq l(m - n) - (n + 1)$. Then,
\begin{align} \begin{split}
    d_1 &\ =\ 2 a_1 + c_1,\\
    d_2 &\ =\ 2 a_2 + c_2.
\end{split} \end{align}

By Lemmas \ref{lemma:loopless orbit length} and \ref{lemma:looped orbit length}, the points in $l(m - n), \dots, m$ that do not constitute the endpoint of a loopless orbit constitute the endpoint of a looped orbit of length $l + 2$, and the points in $n + 1, \dots, l(m - n) - 1$ that do not constitute the endpoint of a loopless orbit constitute the endpoint of a looped orbit of length $l + 1$. Therefore, $c_1$ is the number of looped orbits of length $l + 2$ and $c_2$ is the number of loopless orbits of length $l + 1$.

\begin{proposition} \label{prop:orbit counts}
Depending on the parities of $m, n, l$, the counts of different types of orbits are as follows, and there are no other orbits except the ones listed.
\begin{table}[h!]
\centering
\begin{tabular}{ ||p{0.5cm} p{0.5cm} p{0.5cm}|p{0.5cm} p{0.5cm} p{1.2cm} p{1.2cm}||  }
 \hline
 $m$ & $n$ & $l$ & $c_1$ & $c_2$ & $2 a_1$ & $2 a_2$\\
 \hline \hline
 $1$ & $1$ & $1,0$ & $0$ & $0$ & $d_1$ & $d_2$\\
 $0$ & $0$ & $1,0$ & $1$ & $1$ & $d_1-1$ & $d_2-1$\\
 $1$ & $0$ & $1$ & $1$ & $0$ & $d_1-1$ & $d_2$\\
 $1$ & $0$ & $0$ & $0$ & $1$ & $d_1$ & $d_2-1$\\
 $1$ & $0$ & $1$ & $0$ & $1$ & $d_1$ & $d_2-1$\\
 $1$ & $0$ & $0$ & $1$ & $0$ & $d_1-1$ & $d_2$\\
 \hline
\end{tabular}
\end{table}
\end{proposition}
\begin{proof}
Case 1. (Neither $m$ nor $n$ is even.) By Lemma \ref{lemma:loopless orbit length}, all orbits are loopless, so $c_1 = c_2 = 0$.

Case 2. (Both $m$ and $n$ are even.) The integers $n + 1, \dots, m$ include every number modulo $m - n$ exactly once. Let $q_1$ be the unique integer in $n + 1, \dots, m$ equivalent to $m/2$ modulo $m - n$, and $q_2$ equivalent to $n/2$. Since $m/2 \not\equiv n/2 \pmod {m - n}$, $q_1 \ne q_2$. By Lemma \ref{lemma:looped orbit length}, $q_1$ and $q_2$ constitute the endpoints of looped orbits, and by Lemma \ref{lemma:loopless orbit length}, the orbits that don't end at $q_1$ or $q_2$ are all loopless. Therefore, $c_1 + c_2 = 2$. Since $d_1$ is odd and $2 a_1$ is even, $c_1$ is odd. The only possibility is $c_1 = c_2 = 1$.

Cases 3-6. (One of $m$ and $n$ is even.) If $m$ is even, one orbit has a loop at $m/2$ and the rest are loopless, and if $n$ is even, one orbit has a loop at $n/2$ and the rest are loopless; so $c_1 + c_2 = 1$. When $d_1$ is odd, since $2 a_1$ is even, it has to be that $c_1 = 1$ and $c_2 = 0$. When $d_2$ is odd, since $2 a_2$ is even, it has to be that $c_1 = 0$ and $c_2 = 1$.
\end{proof}

\subsection{Probabilities}
\begin{definition} \label{defn:a_k}
For $k \ge 1$, let $a_k$ be the probability that the chain $(1, \dots, k)$ is satisfied.
\end{definition}

Recalling the definition of a chain, we get the following equivalent definition of $a_k$. If a string of zeros and ones of length $k$ is selected at random such that the probability of a $1$ occurring is $p$, then $a_k$ is the probability that no two consecutive entries are $1$.

\aKRecurrence*

\begin{proof}
The chain $(1)$ is satisfied trivially. The chain $(1, 2)$ is not satisfied if and only if $1$ and $2$ are both in $A$, which happens with probability $p^2$. For $3 \le k \le N$, we note that $(1, \dots, k)$ is satisfied if and only if $(1, \dots, k-1)$ and $(k-1, k)$ are satisfied. With probability $1 - p$, $k$ is non-included in $A$. In that case, $(k-1, k)$ is satisfied with probability $1$ and the probability of satisfaction of $(1, \dots, k-1)$ remains $a_{k-1}$. Therefore, the conditional probability of the satisfaction of $(1, \dots, k)$ given that $k \notin A$ is $a_{k-1}$. With probability $p$, $k$ is included in $A$. In that case, $(k-1, k)$ is satisfied if and only if $k-1$ is non-included in $A$. Therefore, $(1, \dots, k)$ is satisfied if and only if $(1, \dots, k-1)$ is satisfied and $k - 1 \notin A$, which happens if and only if $(1, \dots, k-2)$ is satified and $k - 1 \notin A$. These are independent events that occur with probabilities $a_{k-2}$ and $1 - p$, respectively. Therefore, the conditional probability of the satisfaction of $(1, \dots, k)$ given that $k \in A$ is $(1 - p)a_{k-2}$. Summing over the two possible cases $k \notin A$, $k \in A$, we obtain the desired recurrence. The assumption $k \le N$ is needed because we used $\prob{k \in A} = p$.
\end{proof}

\begin{lemma} \label{lemma:loopless orbit satisfaction prob}
The probability of satisfaction of a loopless orbit $(r_1, \dots, r_k)$ is $a_k$.
\end{lemma}
\begin{proof}
Since $r_1, \dots, r_k$ are distinct, the map $\phi$ that sends $i \in [1, k]$ to $r_i$ and $i \in [k+1, N]$ to $i$ is a bijection from $[1, N]$ to itself. The orbit $(r_1, \dots, r_k)$ is satisfied if for every $i$, at least one of $r_i$ and $r_{i+1}$ is not included in $A$, which is if and only if at least one of $i$ and $i + 1$ is not included in $\phi^{-1}(A)$. Since $\phi^{-1}$ is a bijection, $\phi^{-1}(A)$ has the same probability distribution as $A$.
\end{proof}

\begin{lemma} \label{lemma:looped orbit satisfaction prob}
The probability of satisfaction of a looped orbit $(r_1, \dots, r_{k-1}, r_k = r_{k-1})$ is $(1 - p) a_{k-2}$.
\end{lemma}
\begin{proof}
The chain $(r_1, \dots, r_{k-1}, r_k)$ is satisfied if and only if $(r_1, \dots, r_{k-1})$ and $(r_{k-1}, r_k)$ are both satisfied, which is if and only if $r_{k-1} \notin A$ and $(r_1, \dots, r_{k-1})$ is satisfied, which is if and only if $r_{k-1} \notin A$ and $(r_1, \dots, r_{k-2})$ is satisfied. The two events are independent and occur with probabilities $1 - p$ and $a_{k-2}$, respectively.
\end{proof}

Finally, we can prove Proposition \ref{prop:prob m, n notin sumset}.

\twoNoninclusionProbability*

\begin{proof}
If $[r]_{mn} = (r_1, \dots, r_k)$ and $[q]_{mn} = (q_1, \dots, q_l)$ are two orbits, they are either equal (e.g., $k = l$ and $r_i = q_i$ for all $i$), inverses of each other (e.g., $k = l$ and $r_i = q_{k + 1 - i}$ for all $i$), or disjoint (e.g., $r_i \ne q_j$ for all $i, j$). Since disjoint orbits depend on the inclusion of different numbers into $A$, their satisfactions are independent random events. Orbits that are inverses of each other, $(r_1, \dots, r_k)$ and $(r_k, \dots, r_1)$, are satisfied in all the same circumstances. Theorem \ref{thm:two integers non-included orbits} proves that $m$ and $n$ are non-included in $A + A$ if and only if all orbits are satisfied, and Proposition \ref{prop:orbit counts} counts the different kinds of orbits. The probability that $m$ and $n$ are non-included in $A + A$ is the product of the probabilities of satisfaction of the different disjoint orbits, which are given by Lemmas \ref{lemma:loopless orbit satisfaction prob} and \ref{lemma:looped orbit satisfaction prob}.
\end{proof}

\missingTwoUpperBound*

\begin{proof}
By Lemma \ref{lemma:looped orbit satisfaction prob}, $(1 - p) a_{l - 1}$ is the probability that $(1, \dots, l-1, l-1)$ is satisfied, which is also equal to the probability that $(l+1, l+1, \dots, 2l)$ is satisfied. Therefore, $(1 - p)^2 a_{l-1}^2$ is the probability that $(1, \dots, l-1, l-1, l+1, l+1, \dots, 2l)$ is satisfied, which happens if and only if $(1, \dots, l-1)$, $(l-1, l-1, l+1, l+1)$, and $(l+1, 2l)$ are simultaneously satisfied. By Lemma \ref{lemma:loopless orbit satisfaction prob}, $a_{2l}$ is the probability that $(1, \dots, 2l)$ is satisfied, which happens if and only if $(1, \dots, l-1)$, $(l-1, l, l+1)$, and $(l+1, \dots, 2l)$ are simultaneously satisfied. If $(l-1, l-1, l+1, l+1)$ is satisfied, then $l-1$ and $l+1$ are non-included in $A$, which implies that $(l-1, l, l+1)$ is satisfied. This proves the first inequality. The last two inequalities follow.
\end{proof}

\end{document}